\documentclass[12pt]{amsart}

\newtheorem{theorem}{Theorem}[section]
\newtheorem{proposition}[theorem]{Proposition}
\newtheorem{lemma}[theorem]{Lemma}
\newtheorem{definition}[theorem]{Definition}
\newtheorem{remark}[theorem]{Remark}

\numberwithin{equation}{section}

\begin{document}

\baselineskip=15.5pt

\title{Parabolic principal Higgs bundles}

\author{Indranil Biswas}

\address{School of Mathematics, Tata Institute of Fundamental
Research, Homi Bhabha Road, Bombay 400005, India}

\email{indranil@math.tifr.res.in}



\date{}

\begin{abstract}
In \cite{BBN2}, with Balaji and Nagaraj we introduced the
ramified principal bundles. The aim here is to introduce
the Higgs bundles in the ramified context.
\end{abstract}

\maketitle

\section{Introduction}

Parabolic vector bundles were introduced in \cite{Se}.
The corresponding notion for principal bundles was introduced
in \cite{BBN1} and \cite{BBN2}. In \cite{BBN1}, the parabolic
principal bundles were defined using the Tannakian category
theory. We recall that Nori showed that principal bundles
have a natural description in Tannakian category theory.
The definition of a parabolic principal bundle given in
\cite{BBN1} is modeled on that. In \cite{BBN2} it was
shown that the parabolic principal bundles have a rather
concrete description. More precisely, parabolic principal
bundles were identified with what are called in \cite{BBN2}
as ramified bundles. 

Let $G$ be a complex linear algebraic group.
For a principal $G$--bundle $\psi\, :\, E_G\, \longrightarrow
\, X$, the group $G$ acts freely transitively on each fiber
of $\psi$. For a ramified $G$--bundle $E_G$ over $X$, the
group $G$ acts transitively on each fiber. However, the action
of $G$ on some fibers needs not be free. In other words, points
on some fibers may have nontrivial isotropies;
the details are in \cite{BBN2}, \cite{Bi1}, \cite{Bi2}.
It was shown that there is a natural bijective correspondence
between that parabolic principal $G$--bundle and the
ramified $G$--bundles.

Here we define Higgs fields on a ramified $G$--bundle. For
a principal $G$--bundle $E_G$ over $X$, we recall that a
Higgs field on $E_G$ is a section
$$
\theta\, \in\, H^0(X,\, \text{ad}(E_G)\bigotimes\Omega^1_X)\, ,
$$
where $\text{ad}(E_G)$ is the adjoint vector bundle of $E_G$
and $\Omega^1_X$ is the cotangent bundle, such that
$\theta\bigwedge\theta\, =\, 0$.
We also recall that the adjoint bundle $\text{ad}(E_G)$ is the
one associated to $E_G$ for the adjoint action of $G$ on its
Lie algebra $\mathfrak g$.

The definition of the adjoint vector bundle of a principal
bundle extends to the context of ramified $G$--bundles.
More precisely, for a parabolic $G$--bundle $E_G$
over $X$, its adjoint bundle $\text{ad}(E_G)$
is the parabolic vector bundle over $X$ associated
to $E_G$ for the $G$--module $\mathfrak g$.
However, a Higgs field on a ramified $G$--bundle $E_G$ is not
quite a section of $\text{ad}(E_G)\bigotimes\Omega^1_X$. Rather
it is a section of a bigger sheaf over $X$ containing
$\text{ad}(E_G)\bigotimes\Omega^1_X$, which is called
${\mathcal A}_{E_G}$. Outside the parabolic divisor,
${\mathcal A}_{E_G}$ is identified with
$\text{ad}(E_G)\bigotimes\Omega^1_X$.

\section{Higgs bundle with parabolic structure}\label{sec2}

\subsection{Preliminaries}

Let $X$ be an irreducible smooth projective variety
of dimension $d$ defined over $\mathbb C$.
Let $D\, \subset\, X$ be a simple normal crossing
hypersurface. This means that each irreducible component
of $D$ is smooth and they intersect transversally.
Let $\text{Pvect}(X)$ denote the category
of parabolic vector bundles over $X$ with parabolic structure
over $D$ and rational parabolic weights (see
\cite[Section 2]{Bi2} for more details).

Let $G$ be a complex linear algebraic group.
Let $\text{Rep}(G)$ denote the category of all finite
dimensional rational left representations of $G$. A parabolic
$G$--bundle over $X$ with $D$ as the parabolic divisor is a
functor from $\text{Rep}(G)$ to $\text{Pvect}(X)$ 
that is compatible
with the operations of taking direct sum, tensor product and dual.
(See \cite{BBN1}, \cite[Section 3]{Bi1} and \cite[Section 2]{Bi2};
this approach is based on \cite{No}.)

We will now recall the definition of a ramified $G$--bundles over 
$X$ with ramification over $D$.

A \textit{ramified $G$--bundle} over $X$ with ramification over
$D$ is a smooth complex variety $E_G$ equipped with a right
algebraic action of $G$
$$
f \, :\, E_G\times G\, \longrightarrow\, E_G
$$
and a surjective algebraic map
$$
\psi\, :\, E_G\, \longrightarrow\, X\, ,
$$
such that the following five conditions hold:
\begin{itemize}
\item{} $\psi\circ f \, =\, \psi\circ p_1$, where $p_1$ is
the projection of $E_G\times G$ to $E_G$,

\item{} for each point $x\, \in\, X$, the action of $G$ on the
reduced fiber $\psi^{-1}(x)_{\text{red}}$ is transitive,

\item{} the restriction of $\psi$ to $\psi^{-1}(X\setminus D)$ makes
$\psi^{-1}(X\setminus D)$ a principal $G$--bundle over
$X\setminus D$, meaning the map $\psi$ is smooth over
$\psi^{-1}(X\setminus D)$ and the map to the fiber product
$$
\psi^{-1}(X\setminus D)\times G\, \longrightarrow\,
\psi^{-1}(X\setminus D)\times_{X\setminus D} 
\psi^{-1}(X\setminus D)
$$
defined by $(z\, ,g) \longmapsto\, (z\, , f(z,g))$ is an
isomorphism,

\item{} for each irreducible component $D_i\, \subset\, D$,
the reduced inverse image $\psi^{-1}(D_i)_{\text{red}}$ is a
smooth divisor and
$$
\widehat{D}\, :=\, \sum_{i=1}^\ell \psi^{-1}(D_i)_{\text{red}}
$$
is a normal crossing divisor on $E_G$, and

\item{} for any smooth point $z\in \widehat{D}$, the isotropy
group $G_z\, \subset\, G$, for the action of $G$ on $E_G$,
is a finite cyclic group that acts faithfully on the
quotient line $T_zE_G/T_z\psi^{-1}(D)_{\text{red}}$.
\end{itemize}

There is a natural bijective correspondence
between parabolic $G$--bundles and ramified $G$--bundles.
(See \cite{BBN2}, \cite{Bi2}.) We will interchange without
any further explanation the two terminologies:
ramified $G$--bundle and parabolic $G$--bundle.

The divisor $D$ is fixed once and for all. All parabolic
bundles as well as all the ramified principal bundles are
with respect to
$D$. Hence an explicit mention of $D$ will usually be omitted.

\subsection{Definition of a parabolic principal Higgs bundle}

Let
\begin{equation}\label{psi}
\psi\, :\, E_G\, \longrightarrow\, X
\end{equation}
be a ramified $G$--bundle over $X$ with ramification
over $D$. The Lie algebra of $G$ will be denoted by
$\mathfrak g$. Let
\begin{equation}\label{k0}
{\mathcal K} \, \subset\, TE_G
\end{equation}
be the subbundle defined by the orbits of the action of
$G$ on $E_G$. The action of $G$ on $E_G$ identifies
${\mathcal K}$ with the trivial vector bundle over
$E_G$ with fiber $\mathfrak g$. Let
\begin{equation}\label{eta}
\eta\,:\, E_G\times {\mathfrak g}\, \longrightarrow\,
{\mathcal K}
\end{equation}
be the isomorphism of $\mathcal K$ with the trivial
vector bundle $E_G\times {\mathfrak g}$. This
homomorphism takes the Lie algebra structure of $\mathfrak g$
to the Lie bracket of globally defined vector fields.

Let ${\mathcal Q}$ denote the quotient
vector bundle $TE_G/{\mathcal K}$.
So we have a short exact sequence of vector bundles
\begin{equation}\label{ex.p.0}
0\, \longrightarrow\, {\mathcal K}\, \longrightarrow\, TE_G \,
\stackrel{q}{\longrightarrow}\,{\mathcal Q} \, \longrightarrow\,0
\end{equation}
over $E_G$. The action of $G$ on $E_G$ induces an action of
$G$ on the tangent bundle $TE_G$. This action of $G$ on
$TE_G$ clearly preserves the subbundle ${\mathcal K}$.
It may be mentioned that the isomorphism
$\eta$ in Eq. \eqref{eta} intertwines the action of $G$ on
${\mathcal K}$ and the diagonal action of $G$ constructed
using the adjoint action of $G$ on $\mathfrak g$.
Therefore, we have an induced action of $G$ on the quotient
bundle $\mathcal Q$.

Let
\begin{equation}\label{th0}
\theta_0\, \in\,
H^0(E_G,\, {\mathcal H}om({\mathcal Q}\, ,{\mathcal K}))
\, =\, H^0(E_G,\, {\mathcal K}\bigotimes{\mathcal Q}^*)
\end{equation}
be an algebraic section. We note that the actions of $G$
on ${\mathcal K}$ and ${\mathcal Q}$ together define
an action of $G$ on the complex vector space
$H^0(E_G,\, {\mathcal H}om({\mathcal Q}\, ,{\mathcal K}))$.

Combining
the exterior algebra structure of $\bigwedge
{\mathcal Q}^*$ and the Lie algebra structure on the fibers
of the vector bundle ${\mathcal K}\, =\, E_G\times{\mathfrak g}$
(see Eq. \eqref{eta}), we have a homomorphism
\begin{equation}\label{tau}
\tau\, :\, ({\mathcal K}\bigotimes{\mathcal Q}^*)\bigotimes
({\mathcal K}\bigotimes{\mathcal Q}^*)\, \longrightarrow\,
{\mathcal K}\bigotimes(\bigwedge\nolimits^2{\mathcal Q}^*)\, .
\end{equation}
So $\tau ((A_1\bigotimes\omega_1)\bigotimes(A_2\bigotimes
\omega_2))\,=\, [A_1\, ,A_2]\bigotimes (\omega_1\bigwedge
\omega_2)$. We will denote $\tau (a\, ,b)$ also by $a\bigwedge b$.

\begin{definition}\label{def1}
{\rm A} Higgs field {\rm on a parabolic $G$--bundle
$E_G$ is a section}
$$
\theta_0\, \in\,
H^0(E_G,\, {\mathcal K}\bigotimes {\mathcal Q}^*)
$$
{\rm as in Eq. \eqref{th0} satisfying the following two
conditions:}
\begin{enumerate}
\item {\rm the action of $G$ on $H^0(E_G,\, {\mathcal K}
\bigotimes{\mathcal Q}^*)$ leaves $\theta_0$ invariant, and}
\item $\theta_0\bigwedge\theta_0\, =\, 0$ {\rm (see
Eq. \eqref{tau}).}
\end{enumerate}
\end{definition}

\begin{definition}\label{def1p}
{\rm A} parabolic Higgs $G$--bundle {\rm is a
pair $(E_G\, ,\theta_0)$, where $E_G$ is a
parabolic $G$--bundle, and $\theta_0$ is a
Higgs field on $E_G$.}
\end{definition}

Let
\begin{equation}\label{dag}
{\mathcal A}_{E_G}\, :=\,
(\psi_*({\mathcal K}\bigotimes {\mathcal Q}^*))^G
\end{equation}
be the invariant direct image, where $\psi$ is the
projection in Eq. \eqref{psi}. Therefore,
\begin{equation}\label{e0}
H^0(X,\, {\mathcal A}_{E_G})\, =\,
H^0(E_G,\, {\mathcal K}\bigotimes{\mathcal Q}^*)^G\, .
\end{equation}

For $i\, \geq\, 0$, let
\begin{equation}\label{e00}
\widetilde{\mathcal K}_i\, :=\, (\psi_* ({\mathcal K}
\bigotimes (\bigwedge\nolimits^i{\mathcal Q}^*)))^G
\end{equation}
be the invariant direct image. So,
$\widetilde{\mathcal K}_1\, =\, {\mathcal A}_{E_G}$.
The homomorphism $\tau$
in Eq. \eqref{tau} yields a homomorphism
\begin{equation}\label{e1}
\widetilde{\tau}\, :\, \widetilde{\mathcal K}_1\bigotimes
\widetilde{\mathcal K}_1\, \longrightarrow\,
\widetilde{\mathcal K}_2\, .
\end{equation}

The following lemma is an immediate consequence of the
above constructions.

\begin{lemma}\label{lem1}
A Higgs field on $E_G$ is a section
$$
\theta\, \in\, H^0(X,\, {\mathcal A}_{E_G})
$$
such that
$$
\widetilde{\tau}(\theta\, ,\theta)
\,=\, 0\, ,
$$
where $\widetilde{\tau}$ is constructed in Eq. \eqref{e1}.
\end{lemma}

\subsection{The adjoint vector bundle}

We noted earlier that there is a natural bijective correspondence
between parabolic $G$--bundles and ramified $G$--bundles
(see \cite{BBN2}, \cite{Bi2}). Let $E^P_G$ denote the
parabolic $G$--bundle corresponding to a ramified $G$--bundle
$E_G$. We also recall that $E^P_G$ associates a parabolic
vector bundle over $X$ to each object in $\text{Rep}(G)$. Let
\begin{equation}\label{add0}
\text{ad}(E_G)\, :=\, E^P_G (\mathfrak g)
\end{equation}
be the parabolic vector bundle over $X$ associated to the
parabolic $G$--bundle $E^P_G$ for the adjoint action of
$G$ on its Lie algebra $\mathfrak g$. This parabolic vector
bundle $\text{ad}(E_G)$ will be called the \textit{adjoint
vector bundle} of $E_G$. The vector bundle underlying the
parabolic vector bundle $\text{ad}(E_G)$ will
also be denoted by $\text{ad}(E_G)$. From the
context it will be clear which one is being referred to.

Consider the vector bundle ${\mathcal K}\, \longrightarrow\,
E_G$ constructed in Eq. \eqref{k0}. We noted that ${\mathcal K}$
is equipped with a natural action of $G$.
It is straight forward to check that the
invariant direct image $(\psi_*{\mathcal K})^G$
is identified with the vector bundle underlying the
parabolic vector bundle $\text{ad}(E_G)$ constructed in
Eq. \eqref{add0}.
Indeed, this follows from the fact that for a usual principal
bundle, its adjoint vector bundle coincides with the invariant
direct image of the relative tangent bundle. Therefore, we have
\begin{equation}\label{add}
\text{ad}(E_G)\, =\, (\psi_*{\mathcal K})^G\, .
\end{equation}

There is a natural ${\mathcal O}_X$--linear homomorphism
\begin{equation}\label{a2}
\text{ad}(E_G)\bigotimes \Omega^1_X\, \longrightarrow\,
{\mathcal A}_{E_G}\, ,
\end{equation}
where ${\mathcal A}_{E_G}$ is constructed in Eq. \eqref{dag}.
This homomorphism is an isomorphism over the complement $X
\setminus D$. To prove these, note that
$$
\Omega^1_{X\setminus D}\, =\,
\Omega^1_X\vert_{X\setminus D}\, =\,
(\psi_*{\mathcal Q}^*)^G\vert_{X\setminus D}\, ,
$$
where ${\mathcal Q}$ is the vector bundle in Eq. \eqref{ex.p.0}.
The isomorphism $\Omega^1_{X\setminus D}\, \longrightarrow\,
(\psi_*{\mathcal Q}^*)^G\vert_{X\setminus D}$
extends to a ${\mathcal O}_X$--linear homomorphism
\begin{equation}\label{a22}
\Omega^1_X \, \longrightarrow\,(\psi_*{\mathcal Q}^*)^G
\end{equation}
over $X$. The homomorphism in Eq. \eqref{a2} is obtained from
the isomorphism in Eq. \eqref{add} and the homomorphism
in Eq. \eqref{a22}.

The homomorphism in Eq. \eqref{a2}
is not an isomorphism over $X$ in general.

\section{Semistable parabolic principal Higgs bundles}

\subsection{Reduction of structure group}

Let
\begin{equation}\label{psi2}
\psi\, :\, E_G\, \longrightarrow\, X
\end{equation}
be a ramified principal $G$--bundle with ramification
over $D$. Let
$$
H\, \subset\, G
$$
be a Zariski closed subgroup. Let
$$
U\, \subset\, X
$$
be a Zariski open subset. The inverse image $\psi^{-1}(U)$,
where $\psi$ is the projection in Eq. \eqref{psi2},
will also be denoted by $E_G\vert_U$.

\begin{definition}\label{dr}
{\rm A} reduction of structure group of $E_G$ to $H$
over $U$ {\rm is a subvariety
$$
E_H\, \subset\, E_G\vert_U
$$
satisfying the following three conditions:}
\begin{itemize}
\item {\rm the action of $H$ on $E_G$ preserves $E_H$,}

\item {\rm for each point $x\, \in\, U$, the action of $H$
on $\psi^{-1}(x)\bigcap E_H$ (see Eq. \eqref{psi2})
is transitive, and}

\item {\rm for each point $z\, \in\, E_H$, the isotropy
subgroup of $z$, for the action of $G$ on $E_G$,
is contained in $H$.}
\end{itemize}
\end{definition}

For a point $z\, \in\, E_G$, let
$$
\Gamma_z\, \subset\, G
$$
be the isotropy subgroup of $z$ for the action of $G$ on $E_G$.
It is easy to see that for an element $g\, \in\, G$,
$$
\Gamma_{zg}\, =\, g^{-1}\Gamma_zg\, .
$$
Therefore, for any $z\, \in\, \psi^{-1}(U)$,
if $\Gamma_z\, \subset\, H$, then it follows that
$$
\Gamma_{zg}\, \subset\, H
$$
for all $g\,\in\, H$. Consequently, the last one of the
three conditions in Definition \ref{dr} holds if for each
point $x\, \in\, U$, there exists some point $z\, \in\,
\psi^{-1}(x)\bigcap E_H$ such that
$\Gamma_z\, \subset\, H$.

If
$$
E_H\, \subset\, E_G\vert_U
$$
is a reduction of structure group of $E_G$ to $H$, then
clearly $E_H$ is a ramified principal $H$--bundle over $U$.

\begin{remark}\label{rem1}
{\rm Consider the quotient map}
$$
q_H\, :\, E_G\, \longrightarrow\, E_G/H\, .
$$
{\rm The projection $\psi$ in Eq. \eqref{psi2} defines
a projection}
$$
\psi_H\, :\, E_G/H\, \longrightarrow\, X\, .
$$
{\rm A subvariety}
$$
E_H\, \subset\, E_G\vert_U
$$
{\rm satisfying the first two of the three conditions in
Definition \ref{dr} is given by a section}
$$
\sigma\, :\, U\, \longrightarrow\, E_G/H
$$
{\rm of the projection $\psi_H$. In other words, $\sigma$
is an algebraic morphism and $\psi_H\circ\sigma\, =\,
{\rm Id}_U$. Indeed, the inverse image}
$$
q^{-1}_H(\sigma (U))\, \subset\, E_G
$$
{\rm for any section $\sigma$ satisfies the first two
conditions in Definition \ref{dr}. We note that
such a section $\sigma$ satisfies the third condition
in Definition \ref{dr} if and only if for each point
$x\, \in\, U$, there is a point $z\, \in\,
q^{-1}_H(\sigma (x))$ such that the isotropy subgroup
$\Gamma_z$ of $z$ is contained in $H$. We also note that
if there is one point $z\, \in\,
q^{-1}_H(\sigma (x))$ such that the isotropy subgroup
$\Gamma_z$ is contained in $H$, then the isotropy subgroup
of each point of $q^{-1}_H(\sigma (x))$ is
actually contained in $H$.}
$\hfill{\Box}$
\end{remark}

As before, let $\psi\, :\, E_G\, \longrightarrow\, X$
be a ramified $G$--bundle. Let
\begin{equation}\label{i.}
\iota\, :\, E_H\, \hookrightarrow\, E_G\vert_U
\end{equation}
be a reduction of structure group of $E_G$ to $H$ over a
Zariski open subset $U\, \subset\, X$. Let
\begin{equation}\label{b1}
{\mathcal K}_H\, \subset\, TE_H
\end{equation}
be the subbundle defined by the orbits of the action of $H$
on $E_H$ (see also Eq. \eqref{k0}). Note that
${\mathcal K}_H$ is a subbundle of the pull back $\iota^*
\mathcal K$ (see Eq. \eqref{k0}), where $\iota$ is the
inclusion map in Eq. \eqref{i.}. More precisely,
\begin{equation}\label{i2}
{\mathcal K}_H\, =\, TE_H\bigcap \iota^* \mathcal K\, .
\end{equation}

Let $\mathfrak h$ denote the Lie algebra of $H$.
So $\mathfrak h$ is a submodule of the $H$--module
$\mathfrak g$. The action of $H$ on the Lie algebra
$\mathfrak g$ of $G$ is the restriction of the adjoint
action. The inclusion of the $H$--module $\mathfrak h$ in
$\mathfrak g$ evidently induces an inclusion
\begin{equation}\label{i0}
\text{ad}(E_H)\, \subset\, \text{ad}(E_G)\vert_U
\end{equation}
of parabolic vector bundles (see Eq. \eqref{add0}).

The inclusion of the underlying vector bundles in
Eq. \eqref{i0} also follows from the inclusion of
${\mathcal K}_H$ in $\iota^* \mathcal K$ (see
Eq. \eqref{i2} and Eq. \eqref{add}).

The quotient bundle $(\text{ad}(E_G)\vert_U)/\text{ad}(E_H)$
in Eq. \eqref{i0} has an induced parabolic structure.

On the other hand, consider the ramified principal
$H$--bundle $E_H$ over $U$. Let $E^P_H$ denote the parabolic
principal $H$--bundle over $U$ associated to it. The adjoint
action of the group $H$ on the Lie algebra $\mathfrak g$
produces an action of $H$ on the quotient
${\mathfrak g}/\mathfrak h$,
where $\mathfrak h$ as before is the Lie algebra of $H$. Let
$E^P_H({\mathfrak g}/\mathfrak h)$ be the parabolic
vector bundle over $U$
associated to $E^P_H$ for the $H$--module
${\mathfrak g}/\mathfrak h$. The quotient
parabolic vector bundle $(\text{ad}(E_G)\vert_U)/
\text{ad}(E_H)$ (see Eq. \eqref{i0})
is canonically identified with this
parabolic vector bundle $E^P_H({\mathfrak g}/\mathfrak h)$.
So, we have
\begin{equation}\label{isu}
(\text{ad}(E_G)\vert_U)/
\text{ad}(E_H)\, =\, E^P_H({\mathfrak g}/\mathfrak h)\, .
\end{equation}

Let
$$
{\mathcal Q}_H\, :=\, TE_H/{\mathcal K}_H
$$
be the quotient bundle (see Eq. \eqref{b1}). From Eq.
\eqref{i2} it follows that
${\mathcal Q}_H$ is identified with $\iota^*{\mathcal Q}$, where
${\mathcal Q}$ is the quotient bundle in Eq. \eqref{ex.p.0}.
Therefore, we have a commutative diagram
\begin{equation}\label{mat.}
\begin{matrix}
0 & \longrightarrow & {\mathcal K}_H & \longrightarrow & TE_H &
\longrightarrow & {\mathcal Q}_H & \longrightarrow & 0\\
&& \Big\downarrow && \Big\downarrow && \Vert\\
0 & \longrightarrow & \iota^*{\mathcal K} & \longrightarrow & 
\iota^*TE_G & \longrightarrow & \iota^*{\mathcal Q} &
\longrightarrow & 0
\end{matrix}
\end{equation}
where the bottom exact sequence is pull back of
the one in Eq. \eqref{ex.p.0}
and all the vertical homomorphisms are injective.

Consequently, we have
\begin{equation}\label{i3}
{\mathcal H}om({\mathcal Q}_H\, , {\mathcal K}_H)\, \subset\,
\iota^* {\mathcal H}om({\mathcal Q}\, , {\mathcal K})\, .
\end{equation}

Let
$$
\theta\, \in\,
H^0(E_G,\, {\mathcal H}om({\mathcal Q}\, , {\mathcal K}))
$$
be a Higgs field of $E_G$ (see Definition \ref{def1}).

\begin{definition}\label{de.com.}
{\rm The reduction $E_H$ in Eq. \eqref{i.} is said to be}
compatible {\rm with the Higgs field $\theta$ if}
$$
\theta\vert_{E_H} \, \in\, H^0(E_H,\,
{\mathcal H}om({\mathcal Q}_H\, , {\mathcal K}_H))\,
\subset\, H^0(E_H,\,
\iota^*{\mathcal H}om({\mathcal Q}\, , {\mathcal K}))
$$
{\rm (see Eq. \eqref{i3}).}
\end{definition}

Let
$$
\widehat{\psi}\, :\, E_H\, \longrightarrow\, U
$$
be the natural projection. So, we have $\widehat{\psi}\,=
\,\psi\circ\iota$, where $\iota$ is the inclusion map in Eq.
\eqref{i.} and $\psi$ is the projection of $E_G$ to $X$. Let
\begin{equation}\label{i4}
{\mathcal A}_{E_H} \, :=\, (\widehat{\psi}_*
{\mathcal H}om({\mathcal Q}_H\, , {\mathcal K}_H))^H
\end{equation}
be the invariant direct image, where ${\mathcal Q}_H$ and
${\mathcal K}_H$ are as in Eq. \eqref{mat.} (see also
Eq. \eqref{dag}). From Eq. \eqref{i3} it follows that
\begin{equation}\label{i5}
{\mathcal A}_{E_H}\, \subset\, {\mathcal A}_{E_G}\vert_U\, ,
\end{equation}
where ${\mathcal A}_{E_G}$ is constructed in Eq. \eqref{dag}.

A Higgs field on $E_G$ is an algebraic section
of ${\mathcal A}_{E_G}$ (see Lemma \ref{lem1}). It is now
easy to see that the reduction $E_H$ in Eq. \eqref{i.}
is compatible with a Higgs field
$$
\theta\, \in\, H^0(X, \, {\mathcal A}_{E_G})
$$
on $E_G$ if and only if
$$
\theta\vert_U\, \in\, H^0(U,\, {\mathcal A}_{E_H})\, \subset\,
H^0(U, \, {\mathcal A}_{E_G})
$$
(see Eq. \eqref{i5}).

\subsection{Semistable parabolic principal Higgs bundles}

Fix a very ample line bundle $\xi$ over $X$. The degree
of a torsionfree coherent sheaf $F$ on $X$ is defined to be the
degree of the restriction of $F$ to a smooth complete intersection
curve in $X$ obtained by intersecting hyperplanes on $X$ from the
complete linear system $\vert\xi\vert$. The parabolic
degree of parabolic vector bundle over $X$ is also defined using
$\xi$ (see \cite[p. 81, Definition 1.8(2)]{MY} for the details).
The parabolic degree of a parabolic vector bundle $V$ over $X$
will be denoted by $\mbox{par-deg}(V)$.

In the rest of this section, $G$ will be a
connected reductive linear
algebraic group defined over $\mathbb C$.

Let $(E_G\, ,\theta)$ be a parabolic Higgs $G$--bundle over
$X$. Consider triples of the form $(H\, , U\, , E_H)$, where 
\begin{itemize}
\item $H\, \subset\, G$ is a maximal proper parabolic
subgroup,

\item $U\, \subset \, X$ is a nonempty Zariski open subset
such that the codimension of the complement $X\setminus U$
is at least two, and

\item $E_H\, \subset\, E_G$ is a reduction of structure group
of $E_G$ to $H$ compatible with $\theta$ (see Definition
\ref{de.com.}). 
\end{itemize}

\begin{definition}\label{defst}
{\rm A parabolic Higgs $G$--bundle
$(E_G\, ,\theta)$ over $X$ is called} stable
{\rm (respectively, } semistable{\rm) if for all
triples $(H\, , U\, , E_H)$ of the above type, the inequality}
$$
\mbox{{\rm par}-{\rm deg}}(({\rm ad}(E_G)\vert_U)/
{\rm ad}(E_H))\, >\, 0
$$
{\rm (respectively, $\mbox{{\rm par}-{\rm deg}}(({\rm ad}(E_G)
\vert_U)/{\rm ad}(E_H))\,\geq\, 0$ holds (see Eq. \eqref{isu}).}
\end{definition}

Note that since the codimension of the complement of the open
subset $U$ is at least two, the degree of a vector bundle over
$U$ is well defined.

There is an alternative formulation of the above
definition of (semi)stability which we will now explain.

Let $P$ be a parabolic subgroup of the reductive group $G$.
Therefore, $G/P$ is a complete variety. The quotient map
$G\, \longrightarrow\, G/P$ defines a principal $P$--bundle
over $G/P$. For any character $\lambda$ of $P$, let
$$
L_\lambda\, =\, E_P(\lambda) \, \longrightarrow\, G/P
$$
be the line bundle associated to this principal $P$--bundle
for the character $\lambda$.

Let $Z_0(G)\, \subset\, G$ be the connected component of the
center of $G$ containing the identity element. It is known
that $Z_0(G)\, \subset\, P$. A character $\lambda$ of $P$ which is
trivial on $Z_0(G)$ is called \textit{strictly antidominant}
if the corresponding line bundle $L_\lambda$ over $G/P$ is ample.

Let $(E_G\, ,\theta)$ be a parabolic Higgs $G$--bundle
over $X$. Consider quadruples of the form
$(H\, , \lambda\, ,U\, , E_H)$, where 
\begin{itemize}
\item $H\, \subset\, G$ is a proper parabolic
subgroup,

\item $\lambda$ is a strictly antidominant character of $H$,

\item $U\, \subset \, X$ is a nonempty Zariski open subset
such that the codimension of the complement $X\setminus U$
is at least two, and

\item $E_H\, \subset\, E_G$ is a reduction of structure group
of $E_G$ to $H$ compatible with $\theta$ (see Definition
\ref{de.com.}). 
\end{itemize}

The character $\lambda$ defines an one--dimensional representation
of $H$. Recall that $E_H$ is a parabolic $H$--bundle
over $U$. Let $E_H(\lambda)$ denote the parabolic line bundle
over $U$ associated to the parabolic $H$--bundle
$E_H$ for the $H$--module defined by $\lambda$.

The parabolic Higgs $G$--bundle $(E_G\, ,\theta)$ is
stable (respectively, semistable) if and only if for every
quadruple $(H\, , \lambda\, ,U\, , E_H)$ of the above type,
$$
\text{par-deg}(E_H(\lambda))\, >\, 0
$$
(respectively, $\text{par-deg}(E_H(\lambda))\, \geq\, 0$).

The above assertion follows by imitating the proof of
Lemma 2.1 in \cite[pp. 131--132]{Ra1}.

Let $E_G$ be a parabolic $G$--bundle over $X$.
A reduction of structure group
$$
E_H\, \subset\, E_G
$$
to some parabolic subgroup $H\, \subset\, G$ over $X$ is called
\textit{admissible} if for each character $\lambda$ of
$H$ trivial on $Z_0(G)$, the associated parabolic line
bundle $E_H(\lambda)$ over $X$ satisfies the following
condition:
\begin{equation}\label{admiss}
\text{par-deg}(E_H(\lambda)) \, =\, 0
\end{equation}
(see \cite[p. 307, Definition 3.3]{Ra2} and \cite[pp.
3998--3999]{BS} for admissible reductions of a principal bundle).

A parabolic Higgs $G$--bundle $(E_G\, , \theta)$
over $X$ is called
\textit{polystable} if either $(E_G\, , \theta)$ is stable, or
there is a proper parabolic subgroup $H$ and a reduction of
structure group
$$
E_{L(H)}\, \subset\, E_G
$$
to a Levi subgroup $L(H)$ of $H$ over $X$ such that the following
three conditions hold:
\begin{itemize}
\item the reduction $E_{L(H)}\, \subset\, E_G$ is compatible
with $\theta$,

\item the parabolic Higgs $L(H)$--bundle $(E_{L(H)}\, ,\theta)$
is stable (from the first condition it follows that $\theta$
is also a Higgs field on $E_{L(H)}$), and 

\item the reduction of structure group of $E_G$ to $H$,
obtained by extending the structure group of $E_{L(H)}$
using the inclusion of $L(H)$ in $H$, is admissible.
\end{itemize}

\section{Characteristic classes and connections}

\subsection{Equivariant Higgs $G$--bundles}\label{sec4.1}

Let $Y$ be a complex variety and $\Gamma$ a finite
group acting on $Y$ through algebraic
automorphisms. So we have a homomorphism
\begin{equation}\label{h}
h\, :\, \Gamma\, \longrightarrow\, \text{Aut}(Y)\, ,
\end{equation}
where $\text{Aut}(Y)$ is the group of all automorphisms of
the variety $Y$.

Let $G$ be a complex linear algebraic group.
A $\Gamma$--\textit{linearized} principal $G$--bundle over
$Y$ is a principal $G$--bundle
\begin{equation}\label{ph}
\phi\, :\, F_G\, \longrightarrow\, Y
\end{equation}
and an action of $\Gamma$ on the left of $F_G$
$$
\rho\, :\, \Gamma\times F_G\, \longrightarrow\, F_G
$$
such that the following two conditions hold:
\begin{itemize}
\item the actions of $\Gamma$ and $G$ on $F_G$ commute, and

\item $\phi\circ \rho (\gamma\, ,z)\, =\, h(\gamma)(\phi(z))$
for all $(\gamma\, ,z)\, \in\, \Gamma\times F_G$, where $h$
is the homomorphism in Eq. \eqref{h} and $\phi$ is the projection
in Eq. \eqref{ph}.
\end{itemize}

Let $\psi\, :\, E_G\, \longrightarrow\, X$ be a parabolic
principal Higgs $G$--bundle. There is a finite Galois covering
\begin{equation}\label{vp}
\varphi\, :\, Y\, \longrightarrow\, X
\end{equation}
and a $\Gamma$--linearized principal $G$--bundle $F_G$
over $Y$, where $\Gamma\, :=\, \text{Gal}(\varphi)$ is
the Galois group, such that
\begin{equation}\label{v2}
E_G\, =\, \Gamma\backslash F_G
\end{equation}
(see \cite[Section 3]{Bi2}).

Let $\text{ad}(F_G)$ denote the adjoint vector bundle of
$F_G$. So $\text{ad}(F_G)$ is the vector bundle over $Y$
associated to $F_G$ for the adjoint action of $G$ on
$\mathfrak g$. We recall that an Higgs field on $F_G$
is a section
$$
\theta\, \in\, H^0(Y,\,\text{ad}(F_G)\bigotimes\Omega^1_Y)
$$
such that $\theta\bigwedge\theta\, =\, 0$.

The actions of $\Gamma$ on $Y$ and $F_G$
together induce an action of $\Gamma$ on the vector bundle
$\text{ad}(F_G)\bigotimes\Omega^1_Y$. A Higgs field
$\theta$ on $F_G$ is called $\Gamma$--\textit{invariant} if
the action of $\Gamma$ on $\text{ad}(F_G)\bigotimes\Omega^1_Y$
leaves the section $\theta$ invariant.

\begin{proposition}\label{prop1}
There is a natural isomorphism between the Higgs fields
on $E_G$ and the $\Gamma$--invariant Higgs fields on $F_G$.
\end{proposition}

\begin{proof}
Let
$$
\phi\, :\, F_G\, \longrightarrow\, Y
$$
be the natural projection. Using the action of $G$ on $F_G$,
the kernel of the differential
$$
d\phi\, :\, TF_G\, \longrightarrow\, \phi^*TY
$$
gets identified with the trivial vector bundle $F_G\times
{\mathfrak g}$ over $F_G$ with fiber $\mathfrak g$. Therefore,
\begin{equation}\label{e-1}
H^0(Y,\, \text{ad}(F_G)\bigotimes\Omega^1_Y)\, =\,H^0(F_G,\,
\text{kernel}(d\phi)\bigotimes \phi^*\Omega^1_Y)^G\, .
\end{equation}
Taking the $\Gamma$--invariants of both sides of Eq.
\eqref{e-1} we get a $\mathbb C$--linear isomorphism
$$
H^0(Y,\, \text{ad}(F_G)\bigotimes\Omega^1_Y)^\Gamma\,
\stackrel{\sim}{\longrightarrow}\,
H^0(X,\, {\mathcal A}_{E_G})\, ,
$$
where ${\mathcal A}_{E_G}$ is constructed in Eq. \eqref{dag}.
It is easy to see that the above isomorphism takes the
$\Gamma$--invariant Higgs fields on $F_G$ bijectively
to the Higgs fields on $E_G$.
\end{proof}

Let $(E_G\, ,\theta)$ be a parabolic Higgs $G$--bundle over $X$.
Let $(F_G\, ,\widetilde{\theta})$ be the corresponding
$\Gamma$--linearized principal Higgs $G$--bundle over $Y$ (see
Proposition \ref{prop1}). Fix a polarization $\xi$ on $X$,
and also fix the polarization $\varphi^*\xi$ on $Y$, where
$\varphi$ is the projection in Eq. \eqref{vp}.

\begin{lemma}\label{lem2}
The following three statements are equivalent:
\begin{enumerate}
\item The parabolic Higgs $G$--bundle $(E_G\, ,\theta)$ is
semistable.

\item The $\Gamma$--linearized principal Higgs $G$--bundle $(F_G\, 
,\widetilde{\theta})$ is $\Gamma$--semistable.

\item The principal Higgs $G$--bundle $(F_G\, ,\widetilde{\theta})$
is semistable.
\end{enumerate}
\end{lemma}

\begin{proof}
{}From the second part of Proposition 2.4 of
\cite[p. 26]{BG} it follows that the 
last two statements are equivalent. That the first two
statements are equivalent is obvious.
\end{proof}

\subsection{Characteristic classes}

The adjoint action of $G$ on the Lie algebra $\mathfrak g$ gives
an action of $G$ on ${\mathfrak g}^*$. This in turn defines an
action on the symmetric product $\text{Sym}^n({\mathfrak g}^*)$
for all $n$.

Fix any invariant
\begin{equation}\label{beta}
\beta\, \in\, \text{Sym}^n({\mathfrak g}^*)^G\, .
\end{equation}
For any principal $G$--bundle $F_G$ over $Y$, the
invariant element $\beta$ defines a characteristic class
\begin{equation}\label{beta2}
c_\beta(F_G) \, \in\, H^{2n}(Y,\, {\mathbb C})\, .
\end{equation}
(See \cite[pp. 113--115]{Ch} for the details of the
construction of $c_\beta(F_G)$.)

Now let $Y$ be as in Eq. \eqref{vp}. Let $F_G$ be a
$\Gamma$--linearized principal $G$--bundle on $Y$,
where $\Gamma\, =\, \text{Gal}(\varphi)$. Therefore,
the principal $G$--bundle $h(\gamma)^*F_G$ is isomorphic to
$F_G$ for all $\gamma\, \in\,\Gamma$, where $h$ is the
homomorphism in Eq. \eqref{h}. This implies that
\begin{equation}\label{h2}
h(\gamma)^*c_\beta(F_G)\, =\, c_\beta(F_G)
\end{equation}
for all $\gamma\, \in\,\Gamma$. Consequently, there is a
unique cohomology class
\begin{equation}\label{h3}
\widetilde{c}\, \in\, H^{2n}(Y,\, {\mathbb C})
\end{equation}
such that $\varphi^*\widetilde{c}\, =\, c_\beta(F_G)$,
where $X$ is the quotient $Y/\Gamma$ (see Eq. \eqref{vp}).

Let $E_G\, :=\, \Gamma\backslash F_G$ be the corresponding
parabolic $G$--bundle (see Eq. \eqref{v2}). The
cohomology class $\widetilde{c}$ in Eq. \eqref{h3} will be
called the \textit{characteristic class of $E_G$ for
$\beta$}. The characteristic class of the parabolic 
$G$--bundle $E_G$ for $\beta$ will be denoted by
$c_\beta(E_G)$. So
\begin{equation}\label{h4}
c_\beta(E_G) \, \in\, H^{2n}(X,\, {\mathbb C})\, .
\end{equation}
The integer $n$ in Eq. \eqref{beta} will be called the
\textit{degree} of the characteristic class $c_\beta(E_G)$.

\subsection{Connections and Higgs bundles}

Henceforth, $G$ will be a connected reductive linear algebraic
group defined over $\mathbb C$.

In \cite{Bi2}, we defined holomorphic connections on a
parabolic principal $G$--bundle. A holomorphic connection
on a parabolic principal $G$--bundle is called \textit{flat}
if its curvature vanishes.

\begin{theorem}\label{thm1}
There is a canonical bijective correspondence between
the flat parabolic principal $G$--bundles over $X$ and
the semistable parabolic Higgs $G$--bundles $(E_G\, , \theta)$
over $X$ such that all the characteristic classes of $E_G$
of degree one and degree two vanish.
\end{theorem}

\begin{proof}
Let $E_G$ be a ramified $G$--bundle over $X$. As it was noted
in Section \ref{sec4.1}, there is a finite Galois covering 
$$
\varphi\, :\, Y\, \longrightarrow\, X
$$
and a $\Gamma$--linearized principal $G$--bundle $F_G$
over $Y$, where $\Gamma\, :=\, \text{Gal}(\varphi)$, such that
$$
E_G\, =\, \Gamma\backslash F_G
$$
(see Eq. \eqref{vp} and Eq. \eqref{v2}).

Higgs fields on $E_G$ are the $\Gamma$--invariant Higgs 
fields on $F_G$ (see Proposition \ref{prop1}). Also,
holomorphic connections on the parabolic $G$--bundle $E_G$
are the $\Gamma$--invariant holomorphic connections on $F_G$
(see \cite[p. 269, Proposition 3.4]{Bi1} and
\cite[Theorem 4.4]{Bi2}).

If $\theta$ is a Higgs field on $E_G$ and $\widetilde{\theta}$
the corresponding $\Gamma$--invariant Higgs
fields on $F_G$, then $(E_G\, ,\theta)$ is semistable if and
only if $(F_G\, ,\widetilde{\theta})$ is $\Gamma$--semistable
(see Lemma \ref{lem2}). Also, all the characteristic classes of
$E_G$ of degree one and degree two vanish if and only if
all the characteristic classes of $F_G$
of degree one and degree two vanish.

Therefore, the theorem follows from \cite[p. 20,
Theorem 1.1]{BG}. Note that since $Y$ is a smooth complex
projective variety, $\Gamma$--linearized pseudostable 
Higgs $G$--bundles on $Y$ are the $\Gamma$--linearized
semistable Higgs $G$--bundles on $Y$ (see \cite[p. 26,
Proposition 2.4]{BG}).
\end{proof}

\end{document}